\documentclass{article}
\usepackage{amssymb}
\usepackage{latexsym}
\newcommand{\R}{{\mathbb R}}

\newcommand{\bx}{\hfill{$\Box $ }}

\newtheorem{corollary}{Corollary}[section]
\newtheorem{proposition}{Proposition}[section]
\newtheorem{theorem}{Theorem}[section]

\title{A Detailed Analysis of the Brachistochrone Problem}
\author{R.Coleman\\Laboratoire Jean Kuntzmann,\\ Domaine Universitaire de Saint-Martin-d'H\`eres, France.}
\begin{document}
\date{}
\maketitle

\begin{abstract} The brachistochrone problem gave rise to the calculus of variations. Although its solution is well-known, it is difficult to find a complete and rigourous handling of the problem. The aim of this article is to give a thorough and detailed approach to the brachistochrone problem.  

Classification: 49J05.
\end{abstract}

\hspace{2mm}

If $A$ and $B$ are two points in the plane, with $B$ lower and to the right of
$A$, then we might be tempted to think that an object falling under the
influence of gravity from $A$ would arrive at $B$ most rapidly if it followed
the trajectory of the segment joining $A$ to $B$. Galileo considered this
problem and conjectured that a circular arc would give a better result. Other
scientists over a long period, for example Johann and Jakob Bernouilli, 
Euler and Newton, considered the problem and this eventually gave rise to the
calculus of variations. A solution to this problem is called a brachistochrone. 
In this article we aim to give a rigorous handling of the usual mathematical formulation of this problem. We will assume that all vector spaces are real.\\

We write $C([a,b])$ for the vector space of real-valued continuous functions defined on the closed interval $[a,b]$. The expression  
$$
\|\gamma \| = \sup _{t\in [a,b]}|\gamma (t)|
$$
defines a norm on $C([a,b])$ and with this norm $C([a,b])$ is a Banach space. If  we consider the $y$-axis pointing downwards, then the brachistochrone problem can be formulated in the following way: we take two strictly positive numbers $b$ and $\beta$ and consider the following optimization problem:
$$
\min \frac{1}{\sqrt{2g}}\int _0^b\left(\frac{1+\gamma '^2(t)}{\gamma (t)}\right)^{\frac{1}{2}}dt,
$$
where $g$ is the gravitational constant, $\gamma \in C([0,b])$, $\gamma (0)=0$, $\gamma (b)=\beta$ and $\gamma$ is strictly positive and continuously differentiable on the interval $(0,b]$. This formulation is established in various places, for example \cite{troutman}. In fact, the constant $\frac{1}{\sqrt{2g}}$ plays no role in the search for a minimum, so we can in general neglect it. It should be noticed that the function under the integral sign is not defined at $0$ and so the integral is an improper integral. Hence we need to add the condition that the integral is defined.\\ 

As stated previously, although this problem is well-known, it is difficult to find a complete and rigourous discussion of it. Some or all of the following weaknesses or omissions may be found in most classical texts such as \cite{fox}, \cite{sagan}, \cite{young}:

\begin{itemize}
\item The Euler-Lagrange equation is applied. However, it is overlooked that this is proved for proper integrals not for improper integrals. That it can be adapted to the brachistochrone problem, which involves an improper integral, needs to be proved.
\item It is assumed that an extremal is of class $C^2$. However, this is not necessary, as it can be proved.
\item It is assumed that the Beltrami equation is sufficient to obtain an extremal lying on a cycloid. However, a constant section in a solution is not ruled out by the Beltrami equation. We need to return to the first equation in order to rule out a constant section.
\item In general, very little is said about the cycloid obtained as a function of the second endpoint or where the solution of the brachistochrone problem lies on the cycloid. It is useful to know whether the solution reaches the peak, stops at the peak or goes beyond it. 
\item The extremal obtained is assumed to be the minimum we are looking for. However, this is by no means obvious and needs to be proved. 
\item The relation between the time of transit, i.e., the minimum value of the integral, and the conditions (the pair $(b,\beta )$) is not considered.
\end{itemize}

Given these weaknesses (and possibly others), it seems appropriate to provide a full and rigorous handling of the brachistochrone problem.

\section{Preliminaries} Let $E$ be a vector space and $f$ a real-valued function defined on a nonempty subset $X$ of $E$. Suppose that $v\in E$ and that there exists $\epsilon >0$ such that the segment $[x-\epsilon v, x+\epsilon v]$ is contained in $X$. If the limit
$$
\lim _{t\rightarrow 0}\frac{f(x+tv)-f(x)}{t}
$$
exists, then we call this limit the directional derivative of $f$ at $x$ in the direction $v$ and we write $\partial _vf(x)$ for this limit. The directional derivative is always defined for the vector $0$, but not necessarily for other vectors; however, if it is defined for a certain $v$, then it is also defined for all $\lambda v$ for any $\lambda \in [0,1]$. The directions $v$ for which 
$\partial _vf(x)$ is defined are called ($X$-)admissible directions for $f$ at 
$x$. If $x$ is an extremum (minimum or maximum) of $f$, then $\partial _vf(x)=0$ in all admissible directions for $f$ at $x$. If $E$ is a normed vector space, then this result is also true for local extrema.\\

Let $f$ be a real-valued function defined on a subset $X$ of a vector space $E$ and suppose that, if $x$ and $x+v$ belong to $X$, then the directional
derivative $\partial _vf(x)$ is defined and 
$$
f(x+v)-f(x)\geq \partial _vf(x).
$$
Then we will say that $f$ is convex on $X$. If we have equality only if $v=0$, then we will say that $f$ is strictly convex. Clearly, if $f$ is convex and 
$\partial _vf(x)=0$ for all $v$ such that $x+v\in X$, then $x$ is a minimum, 
which is unique if $f$ is strictly convex.\\

\noindent {\bf Remark.} Usually we define a convex function to be a real-valued 
function $f$ defined on a convex set $X$ such that 
$$
f(x+\lambda v) \leq (1-\lambda )f(x) + \lambda f(x+v)
$$
whenever $x$ and $x+v$ belong to $X$ and $\lambda \in [0,1]$ and we say that $f$ is strictly convex if we have equality only if $v=0$. If $O$ is an open subset of a normed vector space $E$, $X$ a convex subset of $O$ and $f$ a real-valued differentiable function defined on $O$, then $f$ is convex on $X$ if and only 
$$
f(x+v) - f(x) \geq df(x)(v) = \partial _v(x)
$$
whenever $x$ and $x+v\in X$, with strict inequality in the case of strict convexity (see for example \cite{coleman}). This justifies the use of the term convexity above.\\
  
It is often difficult to determine whether a function is convex or not. The
following elementary result, a proof of which may be found in \cite{ciarlet}, 
is very useful in this direction:

\begin{proposition}\label{prop2a} Let $O$ be an open subset of $\R^n$, $X\subset O$ and
$f:O\longrightarrow \R$ of class $C^2$. Then
\begin{itemize}
\item {\bf a.} $f$ is convex on $X$ if and only if the Hessian matrix of $f$ is
positive for all $x\in X$;
\item {\bf b.} $f$ is strictly convex on $X$ if the Hessian matrix of $f$ is
positive definite for all $x\in X$.
\end{itemize}
\end{proposition}

The next elementary result, due to du Bois-Reymond, is fundamental in the 
calculus of variations. Let us write $C^1([a,b])$ for the subspace of 
$C([a,b])$ composed of $C^1$-functions. 

\begin{theorem}\label{theorem1a} If $f\in C([a,b])$ and $\int _a^bf(t)v'(t)dt=0$  for all functions $v\in C^1([a,b])$ such that $v(a)=v(b)=0$, then $f$ is a constant function.
\end{theorem}

\noindent \textsc{proof}  Let $c=\frac{1}{b-a}\int _a^bf(t)dt$ and let us set
$v(s)=\int _a^s(f(t)-c)dt$. Then $v(a)=v(b)=0$ and $v\in C^1([a,b])$,
with $v'(s)=f(s)-c$. Also,
$$
0\leq \int _a^b(f(t)-c)^2dt = \int _a^b(f(t)-c)v'(t)dt
		 = \int _a^bf(t)v'(t)dt-cv(x)|_a^b = 0.
$$
As $f(t)-c$ is continuous, $f(t)-c=0$ for all $t$ and the result follows.\bx

\begin{corollary}\label{cor1a} If $f,g\in C([a,b])$ and 
$$
\int _a^bf(t)v(t) + g(t)v'(t)dt =0
$$
for all functions $v\in C^1([a,b])$ such that $v(a)=v(b)=0$, then $g\in C^1([a,b])$ and $g'=f$.
\end{corollary}

\noindent \textsc{proof} For $s\in [a,b]$ let us set $F(s)=\int _a^sf(t)dt$.
Then $F\in C^1([a,b])$ and $F'(s)=f(s)$. As
$$
\int _a^bf(t)v(t)dt = F(t)v(t)|_a^b - \int _a^bF(t)v'(t)dt = - \int
_a^bF(t)v'(t)dt,
$$
we have
$$
0  = \int _a^bf(t)v(t) + g(t)v'(t)dt = \int _a^b(g(t)-F(t))v'(t)dt.
$$
From the theorem, there is a constant $c\in \R$ such that $g(t)-F(t)=c$, or
$g(t)=F(t)+c$. Therefore $g=F+c\in C^1[a,b]$ and $g'=F'=f$.\bx

\section{Extrema of functionals defined by a definite integral} Suppose that $L$ is a real-valued $C^1$-function defined on an open subset $O\subset\R^2$ and that $\gamma$ is a real-valued, $C^1$-function defined on a closed interval
$\bar{I}=[a,b]$. We also assume that $(\gamma (t),\gamma '(t))\in O$ for all 
$t\in \bar{I}=[a,b]$ and set 
$$
{\cal L}(\gamma ) = \int _a^b L(\gamma (t),\gamma '(t))dt.
$$
To simplify the notation we will write $[\gamma (t)]$ for $(\gamma (t),\gamma '(t))$ and so we may write
$$
{\cal L}(\gamma ) = \int _a^b L[\gamma (t)]dt.
$$
The function $L$ is called a Lagrangian (function). We usually refer to real-valued mappings defined on spaces of functions as functionals. Thus ${\cal L}$ is a (Lagrangian) functional.\\

We now fix $\alpha ,\beta \in \R$ and write $X$ for the subset of $C^1([a,b])$ composed of those $\gamma$ such that $\gamma (a)=\alpha$, $\gamma (b)=\beta$ and $(\gamma (t),\gamma '(t))\in O$ for all $t\in \bar{I}=[a,b]$. These functions form an affine subspace of $C^1([a,b])$. We propose to look for a necessary condition for $\gamma$ to be an extremum of ${\cal L}$ on $X$. To do so, we first find the admissible directions $v$ and an expression for the directional derivative $\partial {\cal L}_v(\gamma )$. If $v$ is an 
admissible direction, then $v\in C^1([a,b])$ and $v(a)=v(b)=0$. In fact, all
such functions $v$ are admissible directions, as we will now see. For $s$ 
small, $\gamma +sv\in X$ and 
$$
\lim _{s\rightarrow 0}\frac{{\cal L}(\gamma +sv)-{\cal L}(\gamma )}{s} = 
\frac{\partial {\cal L}}{\partial s}(\gamma +sv)|_{s=0}.
$$
We have
$$
{\cal L}(\gamma +sv) = \int _a^b L[(\gamma +sv)(t)]dt.
$$
Given the continuity of the integrand with repect to $s$, the derivative 
$\frac{\partial {\cal L}}{\partial s}(\gamma +sv)$ exists for small $s$ 
and so $v$ is an admissible direction. To obtain an expression for the
directional derivative $\partial _v{\cal L}(\gamma )$, we differentiate with 
respect to $s$:
\begin{eqnarray*}
\frac{\partial {\cal L}}{\partial s}(\gamma +sv) &=&
		\int _a^b\frac{\partial L}{\partial s}[(\gamma +sv)(t)]dt\\
		&=& \int _a^b\frac{\partial L}{\partial x}[(\gamma +sv)(t)]v(t)
		+ \frac{\partial L}{\partial y}[(\gamma +sv)(t)]v'(t)dt.
\end{eqnarray*}
As the integrand is continuous with respect to $s$, we obtain
$$
\partial _v{\cal L}(\gamma ) = \int _a^b\frac{\partial L}{\partial x}[\gamma (t)]
v(t)+ \frac{\partial L}{\partial y}[\gamma (t)]v'(t)dt.
$$
Thus we have shown that all $v\in C^1([a,b])$ such that $v(a)=v(b)=0$ are
admissible directions and we have found an expression for the directional
derivative $\partial {\cal L}_v(\gamma )$ for any such $v$. Notice that the admissible directions form a vector subspace of $v\in C^1([a,b])$. \\

If $\gamma$ is an extremum and $v\in C^1[a,b]$ is such that 
$v(a)=v(b)=0$, then $\partial _v{\cal L}(\gamma )=0$, and so from Corollary \ref{cor1a} we obtain
\begin{eqnarray}\label{eqn3a}
\frac{\partial L}{\partial x}[\gamma (t)] &=& 
\frac{d}{dt}\frac{\partial L}{\partial y}[\gamma (t)]
\end{eqnarray}
for $t\in [a,b]$. This equation is known as the Euler-Lagrange equation. Functions which satisfy the Euler-Lagrange equation on some interval are referred to as stationary functions (or extremals). Such functions may or may not be extrema, or even local extrema.\\ 

\section{Extrema of functionals defined by an improper integral} 
In the previous section we supposed that the pair $(\gamma(t),\gamma '(t))$ 
was defined for all $t\in \bar{I}$ and that $(\gamma (t),\gamma '(t))\in O$, the domain of $L$, for all $t\in \bar{I}$. These assumptions are too restrictive to handle the problem which interests us. However, if we slightly relax the conditions, then we still obtain the Euler-Lagrange equation for an extremum.\\ 
 
Let $I$ and $J$ be open intervals of $\R$, where $I=(c,d)$ with $c\in \R$, and $O=I\times J$.  We suppose that $L$ is a real-valued $C^1$-function defined on $O$ and $\gamma\in C([a,b])$ is continuously differentiable on $(a,b]$, with $(\gamma (t),\gamma '(t))\in O$ for all $t\in (a,b]$. If we set 
$$
{\cal L}(\gamma ) = \int _a^b L(\gamma (t),\gamma '(t))dt = \int _a^b L[\gamma
(t)]dt,
$$
then ${\cal L}(\gamma )$ is an improper integral which may or may not be
defined. Let $\alpha ,\beta \in \R$. We will write $X$ for the subset of 
$C([a,b])$ composed of those $\gamma$ which are continuously differentiable on $(a,b]$, with $(\gamma (t),\gamma '(t))\in O$ for all $t\in (a,b]$, and such that $\gamma (a)=\alpha$, $\gamma (b)=\beta$ and 
${\cal L}(\gamma )$ is defined.\\
(Notice that the brachistochrone problem is of this form, with $I=(0,\infty )$, $J=\R$ and $a=\alpha =0$.) 

\begin{theorem} \label{prop3a} If $\gamma $ is an extremum of ${\cal L}$ on 
$X$, then $\gamma$ satisfies the Euler-Lagrange equation on $(a,b]$.
\end{theorem}

\noindent \textsc{proof} Suppose that $v\in C^1([a,b])$ and $v(b)=0$. In addition, assume that there exists $c \in (a,b)$ such that $v$ vanishes on $[a,c]$. Then it is easy to see that $v$ is an admissible direction of ${\cal L}$ at $\gamma$ for any $\gamma\in X$ and 
$$
\partial _v{\cal L}(\gamma) = \int _c ^b\frac{\partial L}{\partial x}[\gamma (t)]
v(t)+ \frac{\partial L}{\partial y}[\gamma (t)]v'(t)dt.
$$
The restriction $u$ of $v$ to $[c,b]$ belongs to $C^1([c,b])$. Therefore, if $\gamma$ is an extremum, then we have
$$
\int _c ^b\frac{\partial L}{\partial x}[\gamma (t)]
u(t)+ \frac{\partial L}{\partial y}[\gamma (t)]u'(t)dt = 0.
$$
We would like to show that this the case for all elements of $C^1([c,b])$ with
$u(c)=u(b)=0$. However, not all members $u$ of $C^1([c,b])$ with 
$u(c)=u(b)=0$ are such restrictions. This will be the case if and only if $u'(c)=0$. Nevertheless, the equality does generally apply as we will now show. 

Let $u\in C^1([c,b])$ with $u(c)=u(b)=0$ and suppose that $u'(c)=\delta >0$. We take $\epsilon \in (0,1]$ such that $d=c-\epsilon >a$ and define a real-valued continuous function $g$ on $[a,c]$ in the following way: $g$ has the value 0 on $[a,d]$, $g$ restricted to $[d,d+\frac{\epsilon}{2}]$ is an inverted `tent' function with height $-\delta$ and $g$ restricted to $[d+\frac{\epsilon}{2},c]$ is an affine function from 0 to $\delta$. If we set
$$
v(t) = 
\cases{\int _a^tg(s)ds & $t\in [a,c]$\cr
u(t) & $t\in [c,b]$,\cr}
$$
then $v$ is a $C^1$-function extending $u$ to $[a,b]$ such that $v$ has the
value 0 on the interval $[a,d]$; hence $v$ is an admissible direction for
${\cal L}$ at $\gamma$. In addition, on the interval $[d,c]$, $|v'(t)|\leq\delta$ and $|v(t)|\leq \frac{\epsilon}{4}\delta \leq\delta$ and so
$$
\left|\int _d ^c\frac{\partial L}{\partial x}[\gamma (t)]
v(t)+ \frac{\partial L}{\partial y}[\gamma (t)]v'(t)dt\right|\leq 
\delta \int _d ^c\left|\frac{\partial L}{\partial x}[\gamma (t)]\right|
+ \left|\frac{\partial L}{\partial y}[\gamma (t)]\right|dt,
$$
which converges to 0, when $\epsilon$ converges to $0$. It now follows that
$$
\int _c ^b\frac{\partial L}{\partial x}[\gamma (t)]
u(t)+ \frac{\partial L}{\partial y}[\gamma (t)]u'(t)dt = 0.
$$
If $u'(c)<0$, then we can use an analogous argument to obtain the same
result. If we now apply Corollary \ref{cor1a}, we see that $\gamma$ satisfies 
the Euler-Lagrange equation on $[c,b]$. As $c$ was chosen arbitrarily in the
interval $(a,b)$, $\gamma$ satisfies the Euler-Lagrange equation on $(a,b]$.\bx\\

The above result gives us a necessary condition for $\gamma$ to be a minimum,
but not a sufficient condition. However, if we add some assumptions, then this 
condition becomes sufficient. Suppose first that $L$ is convex. As $L$ is of
class $C^1$, for any point $x\in O$ the differential $dL(x)$ is defined and
therefore the directional derivative in all directions $h\in \R^2$:
$$
\partial _hL(x) = L'(x)h = \frac{\partial L}{\partial x_1}(x)h_1 + 
\frac{\partial L}{\partial x_2}(x)h_2. 
$$
As $L$ is convex, if $x$ and $x+h$ are in $O$, then
$$
L(x+h) - L(x) \geq \frac{\partial L}{\partial x_1}(x)h_1 + 
\frac{\partial L}{\partial x_2}(x)h_2.
$$
Suppose now that $v$ is of class $C^1$ on $(a,b]$ and such that 
$\int _a^bL[(\gamma+v)](t)dt$ is defined. If $c\in (a,b)$ and $t\in [c,b]$, then
$((\gamma +v)(t),(\gamma +v)'(t))\in O$ and so 
 \begin{eqnarray*}
\int _c^bL[(\gamma +v)(t)]dt - \int _c^bL[\gamma (t)] &\geq& \int
_c^b\frac{\partial L}{\partial x}[\gamma (t)]v(t) + 
\frac{\partial L}{\partial y}[\gamma (t)]v'(t)dt\\
&=& \int _c^b\frac{d}{dt}\frac{\partial L}{\partial y}[\gamma (t)]v(t)+
\frac{\partial L}{\partial y}[\gamma (t)]v'(t)dt\\
&=& \int _c^b\frac{d}{dt}\left(\frac{\partial L}{\partial y}[\gamma (t)]v(t)\right)dt\\
&=& \frac{\partial L}{\partial y}[\gamma (t)]v(t)|_c^b.
\end{eqnarray*}
If we now suppose that $\frac{\partial L}{\partial y}$ is bounded and $\gamma
+v\in X$, then 
$$
{\cal L}(\gamma +v)-{\cal L}(\gamma )=\int _a^bL[(\gamma +v)(t)]dt - \int _a^bL[\gamma (t)] \geq 0,
$$
because $v(a)=v(b)=0$. Therefore $\gamma $ is a minimum. If $L$ is strictly
convex, then an analogous reasoning shows that $\gamma$ is a unique minimum. To
sum up, we have the following result:

\begin{proposition} \label{prop3b} Suppose that $L$ is convex (resp. strictly convex) on $O$
and that $\gamma\in X$ satisfies the Euler-Lagrange equation on $(a,b]$. If 
$\frac{\partial L}{\partial y}$ is bounded, then $\gamma$ is a minimum (resp.
unique minimum) of ${\cal L}$ on $X$.
\end{proposition}

\section{Lagrangians of class $C^2$} We suppose that $L$, $O$ and $X$ are defined as in one of the two previous sections. Our present object is to consider the case where the Lagrangian $L$ is of class $C^2$.  

\begin{theorem} \label{theorem4a} Suppose that $L$ is of class $C^2$ and is 
such that the partial derivative $\frac{\partial ^2L}{\partial y^2}$ does not 
vanish on $O$. If $\gamma \in X$ satisfies the Euler-Lagrange equation, then 
$\gamma$ is of class $C^2$ on $I=(a,b)$.
\end{theorem}

\noindent \textsc{proof} Let us take $t_0\in I$ and set $x_0=\gamma (t_0)$ and 
$y_0=\gamma '(t_0)$. We consider the mapping
$$
\Phi : O \longrightarrow \R \times \R, (x,y) \longmapsto (x, \frac{\partial
L}{\partial y}(x,y)).
$$ 
As $\frac{\partial ^2L}{\partial y^2}(x_0,y_0) \neq 0$, the Jacobian of $\Phi$ 
at $(x_0, y_0)$ is nonzero. From the inverse mapping theorem there is a neighbourhood $U$ of $(x_0,y_0)$ and a neighbourhood $V$ of $(x_0,z_0)$, where $z_0=\frac{\partial L}{\partial y}(x_0,y_0)$, such that $\Phi : U\longrightarrow V$ is a $C^1$-diffeomorphism. We can write
$$
\Phi ^{-1}(x,z) = (x,h(x,z)),
$$
where $h$ is a mapping of class $C^1$. We now define a vector field 
$X:V\longrightarrow \R\times \R$ by
$$
X(x,z) = \left(h(x,z), \frac{\partial L}{\partial x}(x,h(x,z))\right).
$$
$X$ is of class $C^1$, so there is a maximal integral curve $\phi
(t)=(x(t),z(t))$ of $X$, such that $\phi (t_0)=(x_0,z_0)$, defined on an open
interval $J$ containing $t_0$. This integral curve is of class $C^1$. In
addition, $x'(t)=h(x(t),z(t))$ and so $x'(t)$ is of class $C^1$. It
follows that $x(t)$ is of class $C^2$. Let us now set
$$
\psi (t) = (\gamma (t), \frac{\partial L}{\partial y}[\gamma (t)]).
$$
For $t$ close to $t_0$ we have
$$
h\left(\gamma (t),\frac{\partial L}{\partial y}[\gamma (t)]\right) = 
\gamma '(t)
$$
and
$$
\frac{\partial L}{\partial x}\left(\gamma (t), h\left(\gamma (t),\frac{\partial
L}{\partial y}[\gamma
(t)]\right)\right) = \frac{\partial L}{\partial x}[\gamma (t)]=
\frac{d}{dt}\frac{\partial L}{\partial y}[\gamma (t)].
$$
It follows that $\psi$ is an integral curve of $X$. However, $\psi (t_0)=(x_0,z_0)$ and so $\psi (t)= \phi (t)$ on a neighbourhood of $t_0$. Therefore $\gamma (t)=x(t)$ and so $\gamma$ is of class $C^2$ on a neighbourhood of $t_0$. We have shown what we set out to show, namely that $\gamma$ is of class $C^2$ on the interval $I$.\bx\\

Suppose now that $\gamma$ is of class $C^2$, as for example under the 
conditions of the theorem. Then we may derive from the Euler-Lagrange equation
another equation, which is often easier to use. For $t\in (a,b)$ we have

\begin{eqnarray*}
\frac{d}{dt}L[\gamma (t)] &=& \frac{\partial L}{\partial x}
[\gamma (t)]\gamma ' (t) + \frac{\partial L}{\partial y}
[\gamma (t)]\gamma ''(t)\\
 &=& \frac{d}{dt}\frac{\partial L}{\partial y}[\gamma (t)]\gamma '(t)
 + \frac{\partial L}{\partial y}
[\gamma (t)]\gamma ''(t)\\
 &=& \frac{d}{dt}\left(\frac{\partial L}{\partial y}[\gamma (t)]\gamma'(t)\right)
\end{eqnarray*}
and it follows that there is a constant $c$ such that
\begin{eqnarray}\label{eqn4a}
L[\gamma (t)] - \frac{\partial L}{\partial y}[\gamma (t)]\gamma '(t) = c.
\end{eqnarray}
The equation we have just found is called the Beltrami equation.\\

\noindent {\bf Remark.} A function $\gamma$ satisfying the equation $(2)$ is not necessarily a stationary function; however, if $\gamma '$ does not vanish on an
interval, then the Euler-Lagrange equation is satisfied on the interval. Here is
a proof. Suppose that $\gamma '\neq 0$ on an interval $I$ and that $\gamma$
satisfies the equation $(2)$. First, we have
$$
\frac{d}{dt}L[\gamma (t)] = \frac{\partial L}{\partial x}[\gamma (t)]\gamma '(t) + 
\frac{\partial L}{\partial y}[\gamma (t)]\gamma ''(t)
$$
and from equation $(2)$
$$
\frac{d}{dt}L[\gamma (t)] = \frac{d}{dt}\frac{\partial L}{\partial y}[\gamma (t)]\gamma '(t)
 + \frac{\partial L}{\partial y}
[\gamma (t)]\gamma ''(t).
$$
Therefore
$$
\frac{\partial L}{\partial x}[\gamma (t)]\gamma '(t)=
\frac{d}{dt}\frac{\partial L}{\partial y}[\gamma (t)]\gamma '(t)
$$
As $\gamma '(t)\neq 0$, we have
$$
\frac{\partial L}{\partial x}[\gamma (t)] = 
\frac{d}{dt}\frac{\partial L}{\partial y}[\gamma (t)].
$$

\section{The brachistochrone problem and possible solutions} In this section we will apply the previous development to the brachistochrone problem and establish certain properties which a solution must have. For 
$(x,y)\in O=\R _+^*\times \R$ let
$$
L(x,y) = \left(\frac{1+y^2}{x}\right)^{\frac{1}{2}}.
$$
As $$
\frac{\partial L}{\partial x} =
-\frac{1}{2}\left(\frac{1+y^2}{x^3}\right)^{\frac{1}{2}}\qquad\mathrm{and}\qquad 
\frac{\partial L}{\partial y} =
\frac{y}{(x(1+y^2))^{\frac{1}{2}}},
$$
$L$ is of class $C^1$. We fix $b>0$. As in Section 3, for $\gamma\in C([0,b])$ 
continuously differentiable on $(0,b]$, with $(\gamma (t),\gamma '(t))\in O$ for all $t\in (0,b]$, we set 
$$
{\cal L}(\gamma ) = \int _0^b L(\gamma (t),\gamma '(t))dt = \int _0^b L[\gamma
(t)]dt.
$$
The improper integral ${\cal L}(\gamma )$ may or may not be defined. We now take $\beta >0$ and write $X$ for the subset of $C([0,b])$ composed of those $\gamma$ such that ${\cal L}(\gamma )$ is defined, $\gamma (0)=0$ and $\gamma (b)=\beta$. It is easy  to check that, if $\gamma (t)=\frac{\beta}{b}t$, then $\gamma \in X$ and so $X$ is not empty. The brachistochrone problem is to minimize ${\cal L}$ over $X$.\\

There is no difficulty in seeing that the second partial derivatives of $L$ are 
defined and continuous and so $L$ is of class $C^2$. In particular, 
$$
\frac{\partial ^2L}{\partial y^2} = 
	\frac{1}{x^{\frac{1}{2}}(1+y^2)^{\frac{3}{2}}}>0.
$$
As $L$ is of class $C^2$ and $\frac{\partial ^2L}{\partial y^2}\neq 0$, from 
Theorem \ref{theorem4a} a stationary function $\gamma$ is of class $C^2$ and
we may use equation $(2)$. We have
$$
L[\gamma (t)] - \frac{\partial L}{\partial y}[\gamma (t)]{\gamma}'(t) = c,
$$
i.e.,
$$
\left(\frac{1+\gamma'^2(t)}{\gamma (t)}\right)^{\frac{1}{2}} -
\frac{\gamma '^2(t)}{\gamma (t)^{\frac{1}{2}}(1+\gamma'^2(t))^{\frac{1}{2}}} = c,
$$
from which we derive
$$
\frac{1}{\gamma (t)^{\frac{1}{2}}(1+\gamma '^2(t))^{\frac{1}{2}}} = c > 0
$$
and finally
$$
 \gamma (t) (1+\gamma '^2(t)) = k,
$$
where $k=\frac{1}{c^2}$. Any solution of the brachistochrone problem must satisfy such a differential equation on the interval $(0,b)$. Using the Euler-Lagrange equation $(1)$, we can obtain more information. 

\begin{proposition} \label{prop5a} Let $\gamma$ be a solution of the brachistochrone problem.
Then
\begin{itemize}
\item {\bf a.} $\lim _{t\rightarrow 0}\gamma '(t)=\infty$;
\item {\bf b.} $\gamma$ is not constant on an interval;
\item {\bf c.} $\gamma$ has at most one critical point, which is a maximum;
\item {\bf d.} $\gamma$ is either strictly increasing or is unimodal;
\item {\bf e.} $\gamma '$ is strictly decreasing on $(0,b)$. 
\end{itemize} 
\end{proposition}

\noindent \textsc{proof} {\bf a.} It is sufficient to notice that 
$\lim _{t\rightarrow 0}\gamma (t)= 0$.\\

\noindent {\bf b.} From the expression for $\frac{\partial L}{\partial y}$ there  exist continuous functions $a$ and $b$ such that 
$$
\frac{d}{dt}\frac{\partial L}{\partial y}[\gamma (t)]=
\frac{a(t)\gamma ''(t) - \gamma '(t)b(t)}{\gamma
(t)(1+\gamma '^2(t))}.
$$
If $\gamma$ is constant on an interval, then $\frac{d}{dt}\frac{\partial
L}{\partial y}[\gamma (t)]$ vanishes on the interval. However,
the expression $\frac{\partial L}{\partial x}[\gamma (t)]$ does not vanish. It
follows that $\gamma$ is not constant on an interval.\\

\noindent {\bf c.} The function $\gamma$ is bounded by $k$ and reaches the value $k$ at a point $t_0$, if and only if $t_0$ is a critical point. Suppose that  $t_0$ and $t_1$ are both critical points. As $\gamma$ is not constant on the interval $[t_0,t_1]$, there is a point $t$ in the interval such that $\gamma (t)<k$. However, $\gamma$ is continuous on the compact interval $[t_0,t_1]$ and so reaches a minimum at some point $t_2$. As $\gamma (t_2)<k$ and $\gamma '(t_2)=0$, we have a contradiction. Hence there can be at most one critical point, which is clearly a maximum.\\

\noindent {\bf d.} If $\gamma$ has no critical point or has a critical
point at $b$, then $\gamma$ has no critical point in the interval $(0,b)$. 
If there exist points $s$ and $t$, with $s<t$, such that $\gamma (s)=\gamma (t)$, then from Rolle's theorem there exists $r\in (s,t)$, such that $\gamma '(r)=0$, a contradiction. On the other hand, if there exist $s$ and $t$, with $s<t$, such that $\gamma(s)>\gamma (t)$, then from the mean value theorem there exists $v\in (s,t)$ such that $\gamma '(v)<0$. However, as $\gamma (0)=0$ and $\gamma (t)>0$ for all  $t\in (0,b]$, there exists $u\in (0,v)$ such that $\gamma '(u)>0$. From the intermediate value theorem, there exists $r\in (u,v)$, such that $\gamma '(r)=0$, a contradiction. Thus $\gamma $ is strictly increasing.\\
Suppose now that $\gamma$ has a critical point $t'$ in the interval $(0,b)$.
Applying arguments analogous to those which we have just used, we see that $\gamma$ is strictly increasing on the interval $[0,t']$ and strictly decreasing on the interval $[t',b]$.\\

\noindent {\bf e.} This follows directly from {\bf d.} and the differential 
equation satisfied by $\gamma$.\bx\\

\section{Parametric representation of possible solutions} In this section we 
will give a parametric representation of a possible solution $\gamma$ on the
interval $(0,b]$ of the brachistochrone problem and thus learn more about such a possible solution. We set 
$$
h (t) = 
\cases{ 2\arctan\frac{1}{\gamma '(t)} & $t\in (0,b]$ and $\gamma'(t)>0$\cr
\pi & $t\in (0,b]$ and $\gamma '(t)=0$\cr
2(\pi + \arctan\frac{1}{\gamma '(t)}) & $t\in (0,b]$ and $\gamma '(t)<0$\cr}.
$$
(If $\gamma$ does not reach a maximum (resp. reaches a maximum at $b$), then we
ignore the second and third parts (resp. the third part) of the definition.)
It is easy to see that $h$ is continuous and continuously differentiable when $\gamma '$ is nonzero and from Proposition \ref{prop5a} $h$ is strictly 
increasing. It follows that the image of $h$ is an interval $(0,\theta _1]\subset (0,2\pi )$ and $h(b)=\theta _1$. Let us set 
$I_1=(0,\theta _1 )\cap (0,\pi )$ and $I_2=(0,\theta_1 )\cap (\pi ,2\pi )$. 
($I_2$ may be empty.) For $\gamma '(t)\neq 0$ we have
$$
\gamma ' (t) = \cot \frac{h(t)}{2} \Longrightarrow 1+\gamma '^2 (t)=
1+\cot ^2\frac{h(t)}{2} = \frac{1}{\sin ^2\frac{h(t)}{2}}.
$$
Therefore
$$
\gamma (t) = k\sin ^2\frac{h(t)}{2} = \frac{k}{2}(1-\cos h(t)).
$$
Also,
$$
\gamma (t) = \frac{k}{2}(1-\cos h(t))\Longrightarrow \gamma '(t)=
\frac{k}{2}(\sin h(t))h'(t) \Longrightarrow \cot \frac{h(t)}{2} = 
\frac{k}{2}(\sin h(t))h'(t).
$$

Now let us set $h(t)=\theta$. Differentiating $h^{-1}$ on $I_1$ and on $I_2$, if  not empty, we obtain
$$
\frac{d}{d\theta}(h^{-1}) (\theta ) = \frac{k}{2}\frac{\sin\theta}
{\cot\frac{\theta}{2}} = k\sin ^2\frac{\theta}{2} =\frac{k}{2}(1-\cos \theta).
$$
It follows that on the interval $I_1$ (resp. $I_2$, if not empty), 
there is a constant $c_1$ (resp. $c_2$), such that 
$$
t = \frac{k}{2}(\theta -\sin \theta ) + c_i.
$$
As $\lim _{t\rightarrow 0} h(t) = 0$, $c_1=0$ and so on the interval $I_1$
the graph of $\gamma$ has the parametric representation (P):
$$
\cases{t = \frac{k}{2}(\theta  -\sin \theta )\cr
\gamma (t) = \frac{k}{2}(1-\cos \theta ),\cr}
$$
where $\theta \in (0,\pi )$. Suppose that $I_2$ is not empty. As $h^{-1}$ is continuous,
$$
\lim _{\theta \rightarrow \pi +}\frac{h}{2}(\theta -\sin \theta )+c_2 = h^{-1}(\pi ) = \lim _{\theta \rightarrow \pi -}\frac{h}{2}(\theta -\sin \theta ),
$$
therefore $c_2=0$. Thus the parametric representation $(P)$ is valid for the whole graph of $\gamma$: the graph of the function $\gamma$ may be considered as  lying on an arch of a cycloid, i.e., the curve traced out by a point on the cicumference of a disc moving on a plane surface. (The diameter of the disc is $k$.)\\

\noindent {\bf Remark.} As $h(t)=\arccos (1-\frac{2}{k} \gamma (t))$, $h$ is continuously differentiable on its entire domain.\\

The parametric representation we have obtained enables us to learn more about
$\gamma$. Let us consider the function $\alpha $ defined on $(0,2\pi )$ as follows:
$$
\alpha (\theta )= \frac{1-\cos \theta}{\theta -\sin\theta }.  
$$
Clearly $\lim _{\theta \rightarrow 2\pi }\alpha (\theta )=0$. Also 
$$
\theta -\sin\theta = \frac{\theta ^3}{6} + o(\theta ^3) \qquad\mathrm{and}\qquad 1-\cos\theta
=  \frac{\theta ^2}{2} + o(\theta ^3),
$$
therefore $\lim _{\theta \rightarrow 0}\alpha (\theta )=\infty$. A simple
calculation shows that 
$$
\alpha ' (\theta ) = \frac{\theta \sin \theta - 2 + 2\cos \theta}{(\theta -
\sin\theta )^2}.
$$
A careful analysis of the numerator of $\alpha '$ shows that it is strictly 
negative on $(0,2\pi )$ and hence so is $\alpha '$. Therefore $\alpha$ is 
strictly decreasing on $(0,2\pi )$. This means that there is a unique 
$\tilde{\theta}$ such that 
$\frac{\beta}{b}= \alpha (\tilde{\theta })$. However, $\alpha (\theta _1 )=
\frac{\beta}{b}$ and so $\theta _1=\tilde{\theta }$. Therefore, from the value
of $\tilde{\theta }=\alpha ^{-1}(\frac{\beta}{b})$, we may determine whether
$\gamma $ is strictly increasing without a critical point
($\tilde{\theta}<\pi$), strictly increasing with a critical point
($\tilde{\theta}=\pi$) or strictly increasing and then strictly decreasing
($\tilde{\theta}>\pi$). In addition, from one of the equations 
$$
b = \frac{k}{2}(\tilde{\theta}   -\sin \tilde{\theta})\qquad\mathrm{or}\qquad\beta = 
\frac{k}{2}(1-\cos \tilde{\theta}),
$$
we may find $k$. We have shown that if a minimum $\gamma$ of the brachistochrone problem exists, then its graph has a particular form: if 
$$
\tilde{\theta}=\alpha^{-1}(\frac{\beta}{b}) \qquad\mathrm{and}\qquad  
k=\frac{2b}{\tilde{\theta}-\sin \tilde{\theta}},
$$ 
then the graph of $\gamma$ has the parametric representation
$$
\cases{t = \frac{k}{2}(\theta  -\sin \theta )\cr
\gamma (t) = \frac{k}{2}(1-\cos \theta ),\cr}
$$
for $\theta \in [0,\tilde{\theta}]$. We should also notice that the curve so defined is admissible, i.e., it is continuous on $[0,b]$, strictly positive and continuously differentiable on $(a,b]$, and has the endpoint values $0$ and $\beta$. However, we have not shown that it is a minimum. In the next 
section we will look at this question.\\

Before closing this this section we draw attention to two small small details. Firstly, as
$$
\left(\frac{1}{\theta -\sin\theta}\right)' = -\frac{1-\cos\theta}{(\theta -\sin\theta )^2} < 0
$$
for $\theta\in (0,2\pi )$, if $b$ is fixed, then $k$ is an increasing function of $\beta$. (This is of course not surprising.)

Secondly, as $\alpha (\pi )=\frac{2}{\pi}<1$, for the case where $\beta =b$ we must have $\tilde{\theta}<\pi$.\\

\section{Existence of a brachistochrone} From now on we will write $\gamma _0$ for the particuler function we defined parametrically at the end of the last section. We aim to show that $\gamma _0$ is the unique minimum of the brachistochrone problem, i.e., a brachistochrone. We would like to use the criterion developped in Proposition \ref{prop3b}. If 
$$
L(x,y) =  \left(\frac{1+y^2}{x}\right)^{\frac{1}{2}},
$$
then
$$
\frac{\partial ^2L}{\partial x^2}=\frac{3}{4}\frac{(1+y^2)^{\frac{1}{2}}}{x^{\frac{5}{2}}}, \quad  \frac{\partial ^2L}{\partial y^2}=\frac{1}{x^2(1+y^2)^{\frac{3}{2}}} \quad \textrm{and}\quad \frac{\partial ^2L}{\partial y\partial x}=\frac{\partial ^2L}{\partial x\partial y}=-\frac{1}{2}\frac{y}{x^{\frac{3}{2}}(1+y^2)^{\frac{1}{2}}},
$$
and
$$
{\cal H}(L) = \frac{1}{4(1+y^2)}\left(3\frac{1}{x^{\frac{9}{2}}}-\frac{y^2}{x^3}\right),
$$
where ${\cal H}(L)$ is the Hessian of $L$. Clearly ${\cal H}(L)$ is negative for many pairs $(x,y)$ and so by Proposition \ref{prop2a} $L$ is not convex.\\  

We get around this difficulty by introducing another minimization
problem. For $(x,y)\in O=\R _+^*\times \R$, let
$$
M(x,y) = (x^{-2}+y^2)^{\frac{1}{2}}.
$$
As $$
\frac{\partial M}{\partial x} =
-x^{-3}(x^{-2}+y^2)^{-\frac{1}{2}}\qquad\mathrm{and}\qquad 
\frac{\partial M}{\partial y} =y(x^{-2}+y^2)^{-\frac{1}{2}},
$$
$M$ is of class $C^1$. For $\delta\in C[0,b]$ continuously differentiable on 
$(0,b]$, with $(\delta (t),\delta '(t))\in O$ for all $t\in (0,b]$ we set 
$$
{\cal M}(\delta ) = \int _0^b M[\delta (t)]dt.
$$
The improper integral ${\cal M}(\delta )$ may or may not be
defined. We write $Y$ for the subset of $C([0,b])$ composed of those $\delta $ which are continuously differentiable on $(0,b]$ and such that $\delta (0)=0$, 
$\delta (b)=(2\beta )^{\frac{1}{2}}$ and $M(\delta )$ is defined. If $\gamma \in
X$ and we set $\delta =(2\gamma )^{\frac{1}{2}}$, then 
$$
\gamma = \frac{\delta ^2}{2}\qquad\mathrm{and}\qquad \gamma '=\delta \delta '.
$$
It is now easy to check that $\delta\in Y$ if and only if $\delta =(2\gamma )^{\frac{1}{2}}$ for some $\gamma \in X$ and in this case ${\cal L}(\gamma )=2^{\frac{1}{2}}{\cal M}(\delta )$. Let us set $\delta _0=(2\gamma _0)^{\frac{1}{2}}$.

\begin{proposition} $\delta _0$ is the unique minimum of ${\cal M}$ on $Y$.
\end{proposition}

\noindent \textsc{proof} The second partial derivatives of $M$ are defined and 
continuous and so $M$ is of class $C^2$. In fact, 
$$
\frac{\partial ^2M}{\partial x^2}=\frac{2+x^2y^2}{x^6(x^{-2}+y^2)^{\frac{3}{2}}}, \quad  \frac{\partial ^2M}{\partial y^2}=\frac{1}{x^2(x^{-2}+y^2)^{\frac{3}{2}}} \quad \textrm{and}\quad \frac{\partial ^2M}{\partial y\partial x}=\frac{\partial ^2M}{\partial x\partial y}=\frac{y}{x^3(x^{-2}+y^2)^{\frac{3}{2}}},
$$
and
$$
{\cal H}(M) = \frac{2}{x^8(x^{-2}+y^2)^2},
$$
where ${\cal H}(M)$ is the Hessian of $M$. As ${\cal H}(M)$ is positive on $O$, by Proposition \ref{prop2a} $M$ is strictly convex. In addition, $|\frac{\partial M}{\partial y}|<1$. To simplify the notation, let us write $\delta $ for $\delta _0$ and $\gamma$ for $\gamma _0$. We have
$$
\gamma (1+\gamma '^2)=k\Longrightarrow  \frac{\delta ^2}{2}(1 + \delta ^2\delta '^2)=k
$$
and
\begin{eqnarray*}
M[\delta (t)] - \frac{\partial M}{\partial y}[\delta (t)] &=&
(\delta ^{-2}(t) + \delta '^2(t))^{\frac{1}{2}} -
(\delta ^{-2}(t) + \delta '^2(t))^{-\frac{1}{2}}\delta '^2(t)\\
&=& \delta ^{-1}(t)(1+\delta ^2(t)\delta '^2(t))^{\frac{1}{2}}-
\delta (t)(1+\delta ^2(t)\delta '^2(t))^{-\frac{1}{2}}\delta '^2(t)\\
&=& \frac{1}{(2k)^{\frac{1}{2}}}(1+\delta ^2(t)\delta '^2(t))-
\frac{1}{(2k)^{\frac{1}{2}}}\delta ^2(t)\delta '^2(t) \;=\; 
\frac{1}{(2k)^{\frac{1}{2}}}.
\end{eqnarray*} 
From the remark after Theorem \ref{theorem4a} on the interval (resp. two intervals) where $\delta '(t)\neq 0$, $\delta$ satisfies the Euler-Lagrange equation, i.e.,
$$
\frac{\partial M}{\partial x}[\delta (t)] = 
\frac{d}{dt}\frac{\partial M}{\partial y}[\delta (t)].
$$
If $\delta$ has a critical point in the interior of the interval $(0,b)$, then
by continuity the Euler-Lagrange equation is also satisfied at this point.
Therefore the Euler-Lagrange equation is satisfied on $(0,b]$. Applying
Proposition \ref{prop3b} we obtain the result.\bx\\

We are now in a position to show that $\gamma _0$ is the unique solution of the
brachistochrone problem, i.e., a brachistochrone.

\begin{theorem} $\gamma _0$ is the unique minimum of ${\cal L}$ on $X$. 
\end{theorem} 

\noindent \textsc{proof} For $\gamma\in X$, with $\gamma \neq \gamma _0$, we 
have
$$
{\cal L}(\gamma ) = 2^{\frac{1}{2}}{\cal M}((2\gamma )^\frac{1}{2}) > 2^{\frac{1}{2}}{\cal
M}(\delta _0) = {\cal L}(\gamma _0).
$$
This ends the proof.\bx\\

It should be noticed that for distinct pairs $(b_1,\beta _1)$ and 
$(b_2,\beta _2)$ the corresponding solutions $\gamma _1$ and $\gamma _2$ of the
brachistochrone problem are distinct. Let us see why this is so. If $b_1\neq b_2$, then $\gamma _1\neq \gamma _2$, because $\gamma _1$ and $\gamma _2$ are not defined on the same interval. However, it may be so that $\gamma _1$ and $\gamma _2$ lie on the same cycloid. This will be so if $\tilde{\theta}_2=2\pi - \tilde{\theta}_1$ and $\beta _2=\beta _1$, because 
$$
k = \frac{2\beta}{1-\cos\tilde{\theta}}.
$$
Now suppose that $b_1=b_2=b$ and $\beta _1<\beta _2$. If $\gamma _1=\gamma _2=\gamma$, then we have
$$
\frac{2b}{\tilde{\theta}_1-\sin \tilde{\theta}_1} = k = 
\frac{2b}{\tilde{\theta}_2-\sin \tilde{\theta}_2},
$$
which implies that
$$
\tilde{\theta}_1-\sin \tilde{\theta}_1=\tilde{\theta}_2-\sin \tilde{\theta}_2.
$$
As $\theta \longmapsto \theta -\sin \theta$ is an increasing function of $\theta$, this is not possible; therefore, in this case too, $\gamma _1$ and $\gamma _2$ are different.\\

We have just seen that the mapping $F$ from $(\R_+^*)^2$ into $\R_+^*\times (0,2\pi )$ defined by
$$
F(b,\beta ) = (k,\tilde{\theta })
$$
is injective. In fact we can say more.

\begin{theorem} The mapping $F$ is a smooth diffeomorphism.
\end{theorem}

\noindent \textsc{proof} Let $(k,\tilde{\theta})\in \R_+^*\times (0,2\pi )$. If we set 
$$
b = \frac{k}{2}(\tilde{\theta} - \sin\tilde{\theta} )\qquad\textrm{and} \qquad \beta = \frac{k}{2}(1-\cos \tilde{\theta}),
$$
then the pair $(b,\beta )\in (\R_+^*)^2$. As the system of equations 
$$
b = \frac{x}{2}(y - \sin y) \qquad \beta = \frac{k}{2}(1-\cos y)
$$
has a unique solution in $\R_+^*\times (0,2\pi )$, the brachistochrone corresponding to the pair $(b,\beta )$ is defined by by the pair $(k,\tilde{\theta})$ and so the mapping $F$ is surjective. Hence, from what we have seen above, $F$ is bijective.\\

It remains to show that $F$ is smooth. However, it is easy to see that that the inverse mapping $F^{-1}$ has partial derivatives of all orders and so is a smooth mapping. It follows that $F$ also is smooth.\bx\\

\noindent {\bf Remark.} The mapping $F$ gives us a natural identification of brachistochrones with elements of the set $\R_+^*\times (0,2\pi )$.\\

\section{Length of trajectory and time of transit}

To simplify the notation, in this section we write $\gamma$, instead of $\tilde{\gamma}$ for the brachistochrone. We may calculate its length, $l(\gamma )$, using the parametric representation of the curve. We have
\begin{eqnarray*}
l(\gamma ) &=& \frac{k}{2}\int _0^{\tilde{\theta}} \left((1-\cos \theta )^2 + (\sin\theta )^2\right)^{\frac{1}{2}}d\theta \\
	&=& \frac{k}{\sqrt{2}}\int _0^{\tilde{\theta}}(1-\cos \theta )^{\frac{1}{2}}d\theta \\
	&=& k\int _0^{\tilde{\theta}}\sin\frac{\theta}{2}d\theta\\
	&=& 2k(1-\cos\frac{\tilde{\theta}}{2})\\
	&=& 4b\frac{1-\cos\frac{\tilde{\theta}}{2}}{\tilde{\theta}-\sin\tilde{\theta}}.
\end{eqnarray*}
Let us fix $b$ and set 
$$
v(\theta )= \frac{1-\cos\frac{\theta}{2}}{\theta -\sin\theta}.
$$
Then 
$$
v'(\theta ) = \frac{1}{2}\frac{\sin \frac{\theta}{2}}{(\theta -\sin\theta )^2}u(\theta ),
$$
where 
$$
u(\theta ) = \theta + \sin\theta -4\sin\frac{\theta}{2}.
$$
We have
$$
u'(\theta )=1+\cos\theta -2\cos\frac{\theta}{2} = 2\cos\frac{\theta}{2}(\cos\frac{\theta}{2}-1),
$$
and so $u'(\theta )<0$ for $\theta \in (0,\pi )$ and $u'(\theta )>0$ for $\theta \in (\pi ,2\pi )$. Also $u(0)=0$, $u(\pi )=\pi -4$ and $u(2\pi )=2\pi$. It follows that there exists $\theta _1\in (\pi ,2\pi )$ such that $u(\theta )<0$ for $\theta \in (0,\theta _1)$ and $u(\theta )>0$ for $\theta \in (\theta _1,2\pi )$. Thus $v'(\theta )<0$ for $\theta \in (0,\theta _1)$ and $v'(\theta )>0$ for $\theta \in (\theta _1,2\pi )$. If $\alpha (\frac{\beta _1}{b})=\theta _1$, then $l(\gamma)$ is a decreasing function of $\beta$ for $\beta \leq \beta _1$ and an increasing function of $\beta$ for $\beta \geq \beta _1$.\\

The time of transit is given by the value of the integral 
$$
I_{\gamma} = \frac{1}{\sqrt{2g}}\int _0^b\left(\frac{1+\gamma '^2(t)}{\gamma (t)}\right)^{\frac{1}{2}}dt.
$$
We are interested in the minimal time of transit and, in particular, how this varies with the values of the conditions (the pair $(b,\beta )$). Using the notation of the previous section, we have 
$$
1 + {\gamma}'^2(t) = \frac{1}{\sin ^2\frac{\theta}{2}},
$$
with $\theta =\frac{h(t)}{2}$. Also,
$$
\gamma (t) = \frac{k}{2}(1-\cos \theta ) \qquad\mathrm{and}\qquad \frac{dt}{d\theta }=\frac{k}{2}(1-\cos \theta ).
$$
With the variable change $\theta = h(t)$, we obtain
\begin{eqnarray*}
I_{\gamma} &=& \frac{1}{\sqrt{2g}}\int _0^{\tilde{\theta}}\sqrt{k} d\theta\\
	   &=& \frac{1}{\sqrt{2g}}\sqrt{k}\tilde{\theta}\\
	   &=& \frac{1}{\sqrt{2g}}\sqrt{\frac{b}{\gamma}}\frac{\tilde{\theta}}{\sqrt{\tilde{\theta}-\sin\tilde{\theta}}}.
\end{eqnarray*}
Let us fix $b$ and set 
$$
v(\theta )=\frac{\theta}{\sqrt{\theta -\sin\theta}}.
$$
We obtain
$$
v'(\theta ) = \frac{u(\theta )}{(\theta -\sin\theta )^{\frac{3}{2}}}.
$$
where 
$u(\theta )=\frac{\theta}{2}-\sin\theta +\frac{\theta}{2}\cos\theta$. Now, 
$$
u'(\theta )=\frac{1}{2}(1 - \cos\theta -\theta \sin \theta ),
$$
therefore $u'(0)=u'(2\pi )=0$ and $u'(\pi )=1$. Using the fact that 
$u''(\theta )=-\frac{\theta}{2}\cos \theta$, we see that there exists $\theta _2\in (\frac{\pi}{2}, \pi )$ such that $u'$ is negative on $(0,\theta _2)$ and positive on $(\theta _2,2\pi )$. However, $u(0)=u(\pi )=0$ and $u(2\pi )=2\pi$ and so $v'(\theta )<0$ (resp. $>0$) for $\theta \in (0,\pi )$ (resp. $\theta \in (\pi , 2\pi ))$. Therefore, if $\alpha (\frac{\beta _2}{b})=\pi$, then $I_{\gamma}$ is a decreasing function of $\beta$ for $\beta \leq \beta _2$ and an increasing function of $\beta$ for $\beta\geq \beta _2$.\\

It is worth noticing that $\theta _2<\theta _1$, which implies that $\beta _2>\beta _1$. This means that on $(0,\beta _1)$ both $l(\gamma)$ and $I_{\gamma }$ are decreasing functions of $\beta$ and on $(\beta _2,\infty)$ both increasing functions of $\beta$. However, on the interval $(\beta _1, \beta _2)$, $l(\gamma )$ is an increasing function of $\beta$ and $I_{\gamma}$ a decreasing function of $\beta$, i.e., even though the length of the brachistochrone increases, due to its form the transit time decreases.

\hspace{1.25cm}

To conclude, it should be mentioned that recently different approaches to the brachistochrone problem have been developped, for example \cite{balder}, \cite{brookfield}, \cite{lawler}. It is also worth mentioning that the minimization problem has been considered over larger classes of functions, in particular, absolutely continuous functions (for example, see \cite{cesari}). Readable introductions to the history of the brachistochrone problem may be found in \cite {rickey}, \cite{shafer}.  

\hspace{1.25cm}

\noindent {\bf Acknowledgements.} I would like to thank Mohamed El Methni and Sylvain Meignen for their helpful comments and suggestions.

\end{document}